\documentclass
{birkmult}
\usepackage[utf8x]{inputenc}

\usepackage{amsmath}
\usepackage{amssymb}
\usepackage{amsthm}
\usepackage{hyperref}
\usepackage{a4wide}

\newtheorem{thm}{Theorem}[section]
\newtheorem{lem}[thm]{Lemma}

\newtheorem{cor}[thm]{Corollary}

\newtheorem{defi}[thm]{Definition}

\newcommand{\ep}{\varepsilon}
\newcommand{\no}[2]{\|#1\|_{I_{#2},\infty}}
\newcommand{\noo}{\|F_1\|_{I_3,\infty}}

\newcommand{\G}{\ensuremath{{\mathcal G}}}
\newcommand{\NN}{\ensuremath{{\mathcal N}}}
\newcommand{\supp}{\mathrm{supp}}

\newcommand{\Cinf}{\mathcal{C}^\infty}
\newcommand{\D}{\ensuremath{{\mathcal D}}}
\newcommand{\R}{\mathbb{R}}

\newcommand{\grad}{\ensuremath{\mathrm{grad}}}

\newcommand{\E}{\ensuremath{{\mathcal E}}}

\newcommand{\eps}{\varepsilon}
\newcommand{\EM}{\ensuremath{{\mathcal E}_{\mathrm{M}}}}

\newcommand{\comp}{\subset\subset}
\newcommand{\GaG}{\Gamma_\G}
\newcommand{\Ga}{\Gamma}

\renewcommand{\d}{\ensuremath{\partial}}

\begin{document}
\title{Geodesic completeness of generalized space-times}
\author{Clemens S\"amann, Roland Steinbauer}
\address{Faculty of Mathematics, University of Vienna\\Oskar-Morgenstern-Platz 1, 1090 Vienna, Austria}
\email{clemens.saemann@univie.ac.at, roland.steinbauer@univie.ac.at}

\begin{abstract}
We define the notion of geodesic completeness for semi-Riemannian metrics
of low regularity in the framework of the geometric theory of generalized 
functions. We then show completeness of a wide class of impulsive gravitational wave 
space-times.
\end{abstract}

\maketitle

{\bf MSC classes:} 46F30 (Primary), 83C15, 83C35 (Secondary)

\section{Introduction}
The geometric theory of generalized functions (\cite{GKOS01}) based on Colombeau algebras (\cite{C85}) is by now a 
well-established field within generalized functions. It has proved to be widely applicable in geometric situations, 
such as Lie-group analysis of differential equations (e.g.\ \cite{KO00,DKP02,KK06}), wave-type equations on Lorentzian 
manifolds (\cite{GMS09,HKS12}) and various problems in general relativity (see \cite{SV06,NS13} for an overview).

The applications to relativity in particular include the study of the geometry of impulsive gravitational 
waves (introduced by Penrose (\cite{P72}), for a thorough review see \cite{P02}) which are 
key-examples for exact space-times modeling a gravitational wave pulse. The simplest of these
geometries arises as the impulsive limit of  \emph{plane fronted gravitational waves with parallel rays (pp-waves)}
described by the line-element (e.g.\ \cite[Ch.\ 17]{GP09})
\begin{equation}\label{ppw}
 ds^2=-2dudv+dx^2+dy^2+H(x,y,u)du^2,
\end{equation}
on $\R^4$, where $H$ is a smooth function. 
Since the field equations put no restriction on the $u$-dependence of $H$ one can perform the so-called impulsive
limit by basically setting (for details see e.g.\ \cite[Ch.\ 20]{GP09}) 
\begin{equation}\label{ippw}
 H(x,y,u)=f(x,y)\,\delta(u),
\end{equation}
where $f$ is smooth and $\delta$ denotes the 
Dirac-measure. The resulting space-time is flat everywhere but 
on the null-hypersurface $\{u=0\}$, where a gravitational wave impulse is located. 
In \cite{KS99} it was shown that the geodesic equation 
for these geometries possess \emph{unique and globally defined solutions} in nonlinear generalized functions, although 
the global aspect was not emphasized there. It only recently came back into focus 
in the context of causality theory for Lorentzian metrics of low regularity (\cite{CG11,KSS13}).

This work is in particular motivated by a recent result of the authors (\cite{SS12}) which provides a completeness statement 
for a wider class of impulsive radiative geometries, which have been called \emph{impulsive N-fronted waves with parallel rays 
(INPWs)}. These space-times are the impulsive limits of geometries studied originally by Brinkmann in the
context of conformal mappings of Einstein spaces (\cite{B25}) and are of the following form: Let $(N,h)$ be a connected 
Riemannian manifold of dimension $n$, set $M=N\times\mathbb{R}^2_1$ and equip $M$ with the line element
\begin{equation}\label{npw}
 ds^2\ =\ dh^2 + 2dudv + H(x,u)du^2, 
\end{equation}
where $dh^2$ denotes the line element of $(N,h)$. Moreover $(u,v)$ are global
null-coordinates on the $2$-dimensional Minkowski space $\mathbb{R}^2_1$ 
and $H:N\times\mathbb{R}\to\mathbb{R}$ is a smooth function.
The causality and the geodesics of such models have been studied in a series of papers 
(\cite{CFS03,FS03,CFS04,FS06}) since they allow to shed light on some of
the peculiar causal properties of plane waves (i.e., pp-waves
(\ref{ppw}) with $H(x^1,x^2,u)=A_{ij}(u)x^ix^j$), see e.g.\ \cite[Ch.\
13]{BEE96}. 
INPWs now arise as the impulsive limit of \eqref{npw}, i.e., upon setting
\begin{equation}\label{inpwdef}
 H(x,u)=f(x)\delta(u).
\end{equation}

Now for the above mentioned completeness result (\cite{SS12}) the INPW-metric $g$ was replaced by
a net of regularizing metrics $g_\ep$ where the Dirac-$\delta$ is replaced by a strict delta-net
$\delta_\ep$ (for a precise definition see below). More explicitly it deals with the net of 
line elements
\begin{equation}\label{inpw} 
ds_\ep^2\ =\ dh^2 + 2dudv + f(x)\delta_\ep(u)du^2, 
\end{equation} 
on $M$, which physically amounts to viewing the impulsive wave as a limit of
extended sandwich waves with small support but increasing amplitude of the ``profile 
function'' $\delta_\ep$. The result now states that given
any geodesic $\gamma$ in (the smooth) space-time $(M,g_\ep)$ (for $\ep$ fixed) there 
is $\ep_0$ small enough, such that $\gamma$
can be defined for all values of an affine parameter provided $\ep\leq\ep_0$.
Also $\ep_0$, for which the geodesic $\gamma$ becomes complete can be
explicitly estimated in terms of (derivatives of) $f$ and the initial data of
$\gamma$. The obvious drawback of this statement is that $\ep_0$ depends on $\gamma$ i.e., 
in general there is no uniform $\ep_0$ which renders the space-times $(M,g_\ep)$ for fixed $\ep$ geodesically 
complete for all $\ep\leq\ep_0$. 

In this paper we define a notion of completeness for generalized metrics that will allow
to formulate the above completeness result in clear analogy to classical completeness statements. This, in particular,
beautifully exhibits the virtues of a well-founded theory of generalized functions.
\bigskip

For convenience of the reader and to keep this presentation self-contained we start with a 
brief review of semi-Riemannian geometry within the geometric theory of generalized functions.
At the end of section \ref{Sec2} we define geodesic completeness and in section \ref{Sec3} we prove 
geodesic completeness of large classes of impulsive gravitational waves. Finally we discuss associated 
distributions of the global geodesics.

\section{Generalized semi-Riemannian geometry}\label{Sec2}

Colombeau algebras of generalized functions (\cite{C85}) 
are differential algebras which contain the vector space of distributions and display 
maximal consistency with classical analysis. Here we review Lorentzian geometry 
based on the special Colombeau algebra $\G(M)$, for further details see \cite{KS02a,KS02b} and
\cite[Sec.\ 3.2]{GKOS01}. 

Let $M$ be a smooth, second countable Hausdorff manifold.
Denote by $\E(M)$ the set of all nets $(u_\eps)_{\eps\in (0,1]=:I}$ in $\Cinf(M)^I$ 
depending smoothly on $\eps$. Note that smooth dependence on the parameter (which was not assumed in the earlier
references) renders the theory technically more pleasant, while not changing any of the basic properties, see the
discussion in \cite[Section 1]{BK:12}. The \emph{algebra of generalized functions on $M$} (\cite{DD91}) 
is defined as the quotient $\G(M) := \EM(M)/\NN(M)$ of 
\emph{moderate} modulo \emph{negligible} nets in $\E(M)$, where the respective notions are
defined by the following asymptotic estimates 
\[
\EM(M) :=\{ (u_\eps)_\eps\in\E(M):\, \forall K\comp M\
\forall P\in{\mathcal P}\ \exists N:\ \sup\limits_{p\in K}|Pu_\eps(p)|=O(\eps^{-N}) \},\\
\]
\[\NN(M)  :=\{ (u_\eps)_\eps\in\EM(M):\ \forall K\comp M\
\forall m:\ \sup\limits_{p\in K}|u_\eps(p)|=O(\eps^{m}) \},
\]
where ${\mathcal P}$ denotes the space of all linear differential operators on $M$.
Elements of $\G(M)$ are denoted by $u = [(u_\eps)_\eps]$. With componentwise operations and the Lie
derivative with respect to smooth vector fields $\xi\in{\mathfrak X(M)}$ defined by
$L_\xi u :=[(L_\xi u_\eps)_\eps]$, $\G(M)$ is a
\emph{fine sheaf of differential algebras}.

There exist embeddings $\iota$ of $\D'(M)$ into $\G(M)$ that are sheaf homomorphisms
and render $\Cinf(M)$ a subalgebra of $\G(M)$. Another, more coarse way
of relating generalized functions in $\G(M)$ to distributions is based on 
the notion of \emph{association}:
$u\in \G(M)$ is called associated with $v\in \G(M)$, $u\approx v$, if $u_\eps - v_\eps \to 0$
in $\D'(M)$. A distribution $w\in \D'(M)$ is called associated with $u$ if $u\approx \iota(w)$.

The ring of constants in $\G(M)$ is the space $\tilde \R$ of generalized
numbers, which form the natural space of point values
of Colombeau generalized functions. These, in turn, are uniquely
characterized by their values on so-called compactly supported
generalized points.
\smallskip

A similar construction is in fact possible for any locally convex space $F$ in place of $\Cinf(M)$
(\cite{G05}), in particular, $F=\Gamma(M,E)$, the space of 
smooth sections of a vector bundle $E\to M$. The resulting space $\GaG(M,E)$
then is the $\G(M)$-module of \emph{generalized sections}  of the vector bundle
$E$ and can be written as  
\begin{equation}
\label{tensorp} \GaG(M,E) = \G(M) \otimes_{\Cinf(M)} \Ga(M,E)= L_{\Cinf(M)}(\Ga(M,E^*),\G(M)).
\end{equation}
$\GaG$ is a fine sheaf of finitely generated and projective $\G$-modules. 
For the special case of \emph{generalized tensor fields} of rank $r,s$ we use the notation $\G^r_s(M)$, i.e.\ 
\[
 \G^r_s(M)\cong L_{\G(M)}(\G^1_0(M)^s,\G^0_1(M)^r;\G(M)).
\]
Observe that this allows the insertion of generalized vector fields and one-forms into
generalized tensors, which is essential when dealing with generalized
metrics which we define as follows: $g\in\G^0_2(M)$ is called a \emph{generalized pseudo-Riemannian metric}
if it is symmetric ($g(\xi,\eta) = g(\eta,\xi)$ $\forall \xi,\, \eta\in \mathfrak{X}(M)$),
its determinant $\det g$ is invertible in $\G$ (equivalently $|\det (g_\eps)_{ij}| > \eps^m$ 
for some $m$ on compact sets), and it possesses a well-defined 
index $\nu$ (the index of $g_\eps$ equals $\nu$ for $\eps$ small). By a ``globalization Lemma''
in (\cite[Lem.\ 4.3]{HKS12}) any generalized metric $g$ possesses a
representative $(g_\eps)_\eps$ such that each $g_\eps$ is a smooth
metric globally on $M$.

Based on this definition, many notions from (pseudo-)Riemannian geometry can be
extended to the generalized setting. In particular, any 
generalized metric induces an isomorphism between generalized vector fields and
one-forms, and there is a unique Levi-Civita connection $\nabla$ corresponding to $g$.
This provides a convenient framework for non-smooth pseudo-Riemannian geometry and
for the analysis of space-times of low regularity in general relativity which extends
the ``maximal distributional'' setting of \cite{GT87} (\cite{SV09,S08}).
\smallskip

Finally we want to discuss geodesics in generalized semi-Riemannian manifolds
(for details see \cite[Section 5]{KS02b}). To this end we have to introduce \emph{generalized functions taking 
values in the manifold $M$}. More precisely, the space of generalized functions 
defined on a manifold $N$ taking values in $M$, $\G[N,M]$ is again defined as 
a quotient of moderate modulo negligible nets $(f_\ep)_\ep$ of
maps from $N$ to $M$, where we call a net moderate (negligible) if
$(\psi\circ f_\eps)_\eps$ is moderate (negligible) and for all smooth $\psi:M\to\R$.

The \emph{induced covariant derivative} of a generalized vector field $\xi$ on a 
generalized curve $\gamma=[(\gamma_\eps)_\eps]\in\G[J,M]$ (with $J$ a real interval),
can be defined componentwise and gives again a generalized vector field $\xi'$ on
$\gamma$. In particular, a \emph{geodesic} in a generalized pseudo-Riemannian manifold 
is a curve $\gamma\in\G[J,M]$ satisfying $\gamma''=0$. Equivalently the usual
local formula holds, i.e.,  
\begin{equation}\label{geo}
  \Big[\,\Big(\frac{d^2\gamma_\eps^k}{dt^2}
  +\sum_{i,j}\hat\Ga^k_{\eps ij}\frac{\gamma_\eps^i}{dt}\frac{\gamma_\eps^j}{dt}\Big)_\eps\,\Big]
  =0,
\end{equation}
where $\Ga^k_{ij}=[(\Ga^k_{\eps ij})_\eps]$ denotes the Christoffel symbols of the generalized metric
$g=[(g_\eps)_\eps]$.

We now give the following (natural) definition.
\begin{defi}(Geodesic completeness for generalized metrics)
Let $g\in\G_0^2(M)$ be a generalized semi-Riemannian metric. Then the generalized 
space-time $(M,g)$ is said to be \emph{geodesically complete} if every geodesic $\gamma$ 
can be defined on $\R$, i.e., every solution of the geodesic equation
\begin{equation}
 \gamma''=0,
\end{equation}
is in $\G[\R,M]$.
\end{defi}

\section{Geodesic completeness of impulsive gravitational wave space-times}\label{Sec3}

In this section we prove geodesic completeness for a large class of impulsive gravitational waves.
More precisely, we will turn the distributional metrics discussed in the introduction into
generalized metrics and then show that their geodesics can be defined for all values of the
parameter, i.e., they belong to $\G[\R,M]$.

We begin by defining the very general class of regularizations used to turn the distributional
metrics of impulsive wave space-times into generalized metrics.

\begin{defi} A \emph{generalized delta-function} is an element $D\in\G(\R)$ that
has a \emph{strict delta net} $(\delta_\ep)_{\ep\in I}$ as a representative, that
is $(\delta_\ep)_{\ep\in I}$ satisfies the following properties
\begin{center}\label{def-sdn}{
\begin{enumerate}
  \item \label{def-sdn-supp}
    $\supp(\delta_\ep)\subseteq(-\ep,\ep)$ $\forall\ep\in I$,
  \item \label{def-sdn-conv}
    $\int_{\mathbb{R}}{\delta_\ep(x) dx} \to 1$ for $(\ep\searrow 0)$ and
  \item \label{def-sdn-L1b}
  $\exists C > 0: \|\delta_\ep\|_{L^1}=\int_{\mathbb{R}}{|\delta_\ep(x)| dx} \leq C$ $\forall \ep\in I$.
\end{enumerate}}
\end{center}
\end{defi} 

We now may define the generalized metrics used in the following, the
impulsive \emph{pp}-wave (cf.\ \eqref{ppw}, \eqref{ippw})
\begin{equation}\label{gippw}
 ds^2=-2dudv+dx^2+dy^2+f(x,y)D(u)du^2,
\end{equation}
on $M=\R^4$, and the INPW (cf.\ \eqref{npw}, \eqref{inpwdef})
\begin{equation}\label{ginpw}
 ds^2\ =\ dh^2 + 2dudv + f(x)D(u)du^2,
\end{equation}
on $M=N\times\R^2_1$, where as in \eqref{npw}, $(N,h)$ is a connected $n$-dimensional
Riemannian manifold, which from now on we suppose to be \emph{complete}. 
Here $D$ denotes an arbitrary generalized delta function.

For the impulsive pp-waves completeness follows from earlier results. More precisely, by \cite[Thm.\ 1]{KS99}
the geodesic equation has (unique) solutions in $\G[\R,M]$ so that we may state the following result.
\begin{cor}(to \cite[Thm.\ 1]{KS99} --- Completeness of impulsive pp-waves)
 The generalized space-time $(\R^4,g)$ with the metric $g$ given by (\ref{gippw}) is
 geodesically complete.
\end{cor}

This result is of course a special case of completeness of INPWs (just set $N=\R^2$ with the flat metric)
which we prove next. We first have to derive an analog of \cite[Thm.\ 1]{KS99}, which will be
based on \cite[Thm.\ 3.2]{SS12}. 

As detailed in \cite[Sec.\ 2]{SS12} it is possible to choose the coordinate
$u$ as an affine parameter along the geodesics, thereby only
excluding trivial geodesics parallel to the impulse. Hence the geodesic equations of the
space-time \eqref{ginpw} take the form
\begin{eqnarray}\label{eq-geo}
 \ddot{v}(u) &=& -\sum_{j=1}^n\frac{\partial f}{\partial x^j}(x(u))\dot{x}^j(u)\ D(u) -
  \frac{1}{2}\ f(x(u))\ \dot{D}(u),\nonumber\\
 \nabla^{h}_{\dot{x}}\dot{x}(u) &=& \frac{1}{2}\ \grad^h(f(x(u))\ D(u),
\end{eqnarray} 
where $\nabla^h$ denotes the covariant derivative on $(N,h)$, $\grad^h$ is the gradient with 
respect to $h$ on $N$ and $(x^1,\ldots, x^n)$ are coordinates on $N$. We immediately see
that the $v$-equation is linear and decouples from the rest of the system. So it can simply be
integrated and we mainly have to deal with the second equation, which actually is the disturbed geodesic
equation on $N$ with a potential given by $\frac{1}{2}\grad^h(f) D$.

We now give a global existence and uniqueness result for the system \eqref{eq-geo}, where we conveniently 
choose data in front of the impulse at $u=-1$. 

\begin{thm}\label{mainthm}(Existence and uniqueness for geodesics in INPW)
Let $D\in\G(\mathbb{R})$ be a generalized delta function, $f\in\mathcal{C}^\infty(N)$, let $v_0, \dot{v}_0\in\R$, $x_0\in N$ and $\dot{x}_0\in T_{x_0}N$. The initial value problem (\ref{eq-geo}) with data 
\[ v(-1) = v_0,\quad x(-1)=x_0,\qquad \dot{v}(-1)=\dot{v}_0,\quad \dot{x}(-1)=\dot{x}_0,\]
has a unique solution $(v,x)\in \G[\R,\R\times N]$.
\end{thm}
This immediately gives our main result.
\begin{cor}(Completeness of INPW)
 The generalized space-time $(M,g)$ with the metric $g$ given by (\ref{ginpw}) is
 geodesically complete.
\end{cor}

In the (uniqueness part of the) proof of the theorem we need the following variant of \cite[Lem.\ A.2]{SS12}. 
\begin{lem}\label{Lem-Basic-Ex}
Let $F_1\in \mathcal{C}^{\infty}(\mathbb{R}^{2n},\mathbb{R}^n)$, $F_2 \in \mathcal{C}^{\infty}(\mathbb{R}^n,\mathbb{R}^n)$,
let $J:=[-1,1]$, $k\in\Cinf(J, \R^n)$ be bounded, let $x_0,\dot{x}_0\in\mathbb{R}^n$, let $b>0$, $c>0$ be given and let
$(\delta_\ep)_\ep$ be a strict delta net with $L^1$-bound $C>0$. Define $I_1:=\{x\in\mathbb{R}^n: |x-x_0|\leq b \}$,
$I_2:=\{x\in\mathbb{R}^n:|x-\dot{x}_0|\leq c + C \no{F_2}{1} \}$ and $I_3:=I_1\times I_2$. Moreover set 
\begin{equation}
\alpha:=\min(1,\frac{b}{|\dot{x}_0| + \noo+C \no{F_2}{1} + \|k\|_{J,\infty}},\frac{c}{\noo + \|k\|_{J,\infty}})\ .
\end{equation}
Then the regularized problem
\begin{equation}\label{eq-geo1}
\left\{ \begin{array}{ll}
\ddot{x}=F_1(x,\dot{x})+F_2(x)\delta_\epsilon + k, \\
x(-\epsilon)=x_0,\ \dot{x}(-\epsilon)=\dot{x}_0,
\end{array} \right.
\end{equation}
has a unique solution $x_\ep$ on $J_{\epsilon}:=[-\epsilon,\alpha-\epsilon]$. Moreover $x_\ep$ and
$\dot x_\ep$ are locally uniformly bounded. Finally the result remains true if we replace $k$ by
a net $(k_\ep)_\ep$ in $\Cinf(J,\R^n)$ which is uniformly bounded.
\end{lem}

The proof is obtained by adapting the proof of \cite[Lem.\ A.2]{SS12}, so that it is not 
necessary to give it here. Note that by classical ODE-theory Lemma \ref{Lem-Basic-Ex} gives global uniqueness (not
only in the function space $X_\ep$ used in the proof of \cite[Lem.\ A.2]{SS12}). 
\medskip

{\bf Proof of the Theorem.}
We proceed in three steps. \emph{First} we need to obtain a "solution candidate", i.e., a net of 
smooth solutions $(v_\eps, x_\eps)_\eps$ defined for all of $\R$  (at least for small $\eps$)
of the regularized initial value problem
\begin{eqnarray}\label{eq-geo-reg}
 \ddot{{v}}_\ep &=& -\delta_\ep\sum_{j=1}^n \frac{\d f}{\d{x^j}}({x}_\ep)\dot{{x}}_\ep^j
-\frac{1}{2}f({x}_\ep)\dot{\delta}_\ep,\nonumber \\
\ddot{{x}}_\ep^k &=& -\sum_{i,j=1}^n\Gamma_{ij}^{k(N)}({x}_\ep)\dot{{x}}_\ep^i
\dot{{x}}_\ep^j + \frac{1}{2} \delta_\ep\sum_{m=1}^n h^{km}({x}_\ep)\frac{\d f}{\d{x^m}}({x}_\ep),\\
{v}_\ep(-1) &=& v_0, \quad {x}_\ep(-1)\,=\,x_0,\quad
\dot{{v}}_\ep(-1)=\dot{v}_0, \quad \dot{{x}}_\ep(-1)\,=\,\dot{x}_0,\nonumber
\end{eqnarray}
where $\Gamma_{ij}^{k(N)}$ denotes the Christoffel symbols of the (smooth) metric $h$ on $N$.
By \cite[Thm.\ 3.2]{SS12} we obtain such a solution for all $\eps$ smaller than a certain 
$\eps_0$ (which depends on $\dot x_0$ as well as on $\Gamma_{ij}^{k(N)}$ and $\grad^h(f)$
on a neighborhood of the point where the geodesic crosses the impulse). Since we are interested
in generalized solutions we may choose $(v_\eps,u_\eps)$ arbitrarily (yet smoothly depending on $\ep$) for all $\eps_0\leq
\eps\leq 1$.
 
\medskip

In the \emph{second} step we prove \emph{existence} of solutions in $\G[\R,\R\times N]$, that
is we establish that the net obtained in the first step is moderate. To this end we have to show
that $v_\eps$ and  $\psi\circ x_\eps$ are moderate for arbitrary smooth $\psi\colon N\to\R$. 
But by Lemma \ref{Lem-Basic-Ex} (with with $b>0$, $c>0$, 
$k=0$, $F_1(y,z)^k:=-{\Gamma_{ij}^{k(N)}(y)z^i z^j}$, $F_2(y)^k:=\frac{1}{2} {h^{km}(y)\frac{\d f}{\d{x^m}}(y)}$)
the solution $x_\eps$ and its derivative $\dot x_\eps$ are even uniformly bounded on compact subsets of $\R$. Using the
differential equation inductively, we see that in fact the derivative of order $l$ of $x_\eps$ is of order $O(\eps^{1-l})$.
Now from the $v$-equation it follows that $\ddot v_\eps$ obeys an $O(\eps^{-2})$-estimate
and inductively all higher order derivatives obey $O(\eps^{-l})$-estimates. The estimates for
$v_\eps$ and $\dot v_\eps$ simply follow by integration. Smooth dependence on $\ep$ is immediate.
\medskip

In the \emph{third} step it remains to show \emph{uniqueness}. To this end 
suppose that $(\tilde{v},\tilde{x})\in\G[\R,\R\times N]$ is a solution of \eqref{eq-geo} as well. 
Writing $\tilde{v}=[(\tilde{v}_\ep)_\ep]$ and $\tilde{x}=[(\tilde{x}_\ep)_\ep]$, 
there exist $[(a_\ep)_\ep]\in \NN(\R)$, $[(b_\ep)_\ep]\in\NN(N)$ and negligible
generalized numbers $[(c_\eps)_\eps]$, $[(\dot c_\eps)_\eps]$, $[(d_\eps)_\eps]$,  
$[(\dot d_\eps)_\eps]$, such that
\begin{eqnarray*}\label{eq-geo-ep}
\ddot{\tilde{v}}_\ep &=& -\delta_\ep\sum_{j=1}^n{ \frac{\d f}{\d{x^j}}(\tilde{x}_\ep)\dot{\tilde{x}}_\ep^j}
-\frac{1}{2}f(\tilde{x}_\ep)\dot{\delta}_\ep + a_\ep,\\
\ddot{\tilde{x}}_\ep^k &=& -\sum_{i,j=1}^n{\Gamma_{ij}^{k(N)}(\tilde{x}_\ep)\dot{\tilde{x}}_\ep^i
\dot{\tilde{x}}_\ep^j} + \frac{1}{2} \delta_\ep\sum_{m=1}^n{h^{km}(\tilde{x}_\ep)\frac{\d f}{\d{x^m}}(\tilde{x}_\ep)} +
b_\ep^k,\\
\tilde{v}_\ep(-1) &=& v_0 + c_\ep, \quad \tilde{x}_\ep(-1)\, =\, x_0 + d_\ep,\quad
\dot{\tilde{v}}_\ep(-1)\,=\,\dot{v}_0 + \dot{c}_\ep, \quad \dot{\tilde{x}}_\ep(-1)\,=\,\dot{x}_0 + \dot{d}_\ep.
\end{eqnarray*}

We have to show that $(x_\eps-\tilde x_\eps)_\eps$ is negligible. To this end we will, however, also estimate
$(\dot x_\eps-\dot{\tilde{x}}_\eps)_\eps$. 

Now by Lemma \ref{Lem-Basic-Ex} we know that $(x_\ep)_\ep$, $(\tilde{x}_\ep)_\ep$ are locally uniformly
bounded. Moreover using the integral
formulas for $x_\ep$ respectively for $\tilde{x}_\ep$ and that $(a_\ep)_\ep$, $(b_\ep)_\ep$ are negligible, 
we obtain that

\begin{align*}
\forall T>0\ &\forall q\in\mathbb{N}\ \exists K_1, K_2>0\ \exists \eta>0 \text{ such that } \forall \ep\in(0,\eta)\
\forall u\in[-T,T]:\\
&|x_\ep^i(u)-\tilde{x}_\ep^i(u)|\leq K_1 \ep^q + \int_{-\ep}^u\int_{-\ep}^s |F_1(x_\ep(r),\dot{x}_\ep(r)) -
F_1(\tilde{x}_\ep(r),\dot{\tilde{x}}_\ep(r))| dr ds\ +\\
&\qquad\int_{-\ep}^u\int_{-\ep}^s |F_2(x_\ep(r)) - F_2(\tilde{x}_\ep(r))| |\delta_\ep(r)| dr ds \leq\\
& K_1 \ep^q + C_3 \int_{-\ep}^u\int_{-\ep}^s(|x_\ep(r)-\tilde{x}_\ep(r)| + |\dot{x}_\ep(r)-\dot{\tilde{x}}_\ep(r)|) dr ds\
+\\
&\qquad C_4 \int_{-\ep}^u\int_{-\ep}^s|x_\ep(r)-\tilde{x}_\ep(r)| |\delta_\ep(r)| dr ds.\\
&\hspace{-22mm}\text{Similarly }\text{for the derivatives: }\\
&|\dot{x}_\ep^i(u)-\dot{\tilde{x}}_\ep^i(u)|\leq K_2 \ep^q + \int_{-\ep}^u|F_1(x_\ep(s),\dot{x}_\ep(s)) -
F_1(\tilde{x}_\ep(s),\dot{\tilde{x}}_\ep(s))| ds\ +\\
&\qquad\int_{-\ep}^u |F_2(x_\ep(s)) - F_2(\tilde{x}_\ep(s))| |\delta_\ep(s)| ds \leq\\
& K_2 \ep^q + C_3 \int_{-\ep}^u(|x_\ep(s)-\tilde{x}_\ep(s)| + |\dot{x}_\ep(s)-\dot{\tilde{x}}_\ep(s)|) ds\ +\\
&\qquad C_4 \int_{-\ep}^u(|x_\ep(s)-\tilde{x}_\ep(s)| |\delta_\ep(r)| ds.
\end{align*}
Here we have used the mean value theorem to obtain the constants $C_3, C_4$. Adding these two inequalities and setting
$\psi(u):=|x_\ep(u)-\tilde{x}_\ep(u)| + |\dot{x}_\ep(u)-\dot{\tilde{x}}_\ep(u)|$ (for $u\in [-T,T]$) yields
\begin{align*}
 \psi(u) \leq (K_1+K_2)\ep^q +
\int_{-\ep}^u(C_3 + C_4|\delta_\ep(s)|)\psi(s) ds + \int_{-\ep}^u\int_{-\ep}^s(C_3+C_4 |\delta_\ep(r)|)\psi(r) dr ds\ .
\end{align*}
Then by a generalization of Gronwall's inequality (due to Bykov \cite[Thm.11.1]{BS:92}) we get that
\begin{align*}
\psi(u)\leq (K_1+K_2)\ep^q \exp\left(\int_{-\ep}^u(C_3 + C_4 |\delta_\ep(s)|)ds + \int_{-\ep}^u\int_{-\ep}^s(C_3 + C_4
|\delta_\ep(r)|)dr ds\right)\leq K'\ep^q,
\end{align*}
where we used the fact that $u\in[-T,T]$ and the uniform $L^1$-bound on $(\delta_\ep)_\ep$, i.e.,
\ref{def-sdn},\ref{def-sdn-L1b}.
This shows that $(x_\ep-\tilde{x}_\ep)$ is negligible (by \cite[Thm. 1.2.3]{GKOS01}). Furthermore since
$(v_\ep-\tilde{v}_\ep)$ can be obtained by integrating $(x_\ep-\tilde{x}_\ep)$ we conclude that it is negligible too.
\hfill $\Box$

\bigskip

Finally we briefly discuss what can be said classically about the geodesics,
that is we provide associated distributions for the unique global geodesics
$(v,x)\in\G[\R,\R\times N]$ obtained in Theorem \ref{mainthm}. In \cite[Sec.\ 4]{SS12} it was shown that
\[
  x_\ep\,\approx\,y, \quad
  v_\ep\,\approx\, w -\frac{1}{2}\ f(x(0))H -\sum_{j=1}^n\Big(\dot{x}^j(0)+\frac{1}{4}\grad^h(f)^j(x(0))\Big) \d_j f(x(0))u_+,
\]
where the first relation even holds in the sense of $0$-association, i.e., the convergence is locally
uniformly. Here the limit $y$ is given by pasting together appropriate (unperturbed) geodesics of the
background $(N,h)$, i.e.,
\begin{equation}\label{ro:xlim}
 y(u):=\left\{\begin{array}{ll}x(u)&u\leq 0,\\\tilde x(u)&u\geq 0,\end{array}\right.
\end{equation}
where $x$ and $\tilde x$ are solution of $\nabla^h_{\dot x}\dot x=0$ with data $x(-1)=x_0$, 
$\dot{x}(-1)=\dot{x}_0$, and $\tilde{x}(0)=x(0)$, $\dot{\tilde{x}}(0)=\dot x(0)+ \frac{1}{2}\grad^h(f)(x(0))$ 
respectively. Moreover $w(u) = v_0 + \dot{v}_0(1 + u)$, $H$ is the Heaviside function and $u_{+}(u)=uH(u)$ 
denotes the kink function.

Hence the $x$-component is continuous, while in general it is not 
differentiable at the point where it hits the impulse, i.e., at $u=0$, with its 
derivative having a jump there. Furthermore the $v$-component is discontinuous, 
it has a  jump at the impulse and a $\delta$ in its derivative. 
We also observe that the parameters of the jump and the refraction at $u=0$ are given 
in terms of $f$ and its derivative at the point where the geodesic hits the impulse.

\subsubsection*{Acknowledgment}
This work was supported by projects P23714, P25326 of the Austrian Science Fund and OEAD WTZ Project CZ15/2013.

\bibliographystyle{alpha}
\bibliography{iNPWsC}
\addcontentsline{toc}{section}{References}

\end{document}